\documentclass[12pt, a4paper,twosided,reqno]{amsart}
\usepackage{enumerate, cite}
\usepackage{hyperref}

\newtheorem{theorem}{Theorem}
\newtheorem{lemma}{Lemma}

\newtheorem{example}{Example}
\newtheorem{corollary}{Corollary}

\newtheorem{remark}{Remark}

\allowdisplaybreaks
\author[G. Pant, S.K.Pant]{Garima Pant, Sanjay Kumar Pant}
\address{Garima Pant; department of mathematics, university of delhi, delhi-110007, india.}
\email{garimapant.m@gmail.com}
\address{Sanjay Kumar Pant; department of mathematics, deen dayal upadhyaya college, university of delhi, new delhi-110078, india.}
\email{skpant@ddu.du.ac.in}
\thanks {Research work of the first author is supported by research fellowship from University Grants Commission (UGC), New Delhi, India.}
\title[entire solution of certain...]{ Entire Solutions of Certain Non-Linear Differential-Difference Equations}
\subjclass{ 34M05, 30D35, 39B32}
\keywords {Nevanlinna theory, Entire function, Difference equation,  differential-difference equation}
\begin{document}
\maketitle	
\begin{abstract}
In this paper, we study about existence and non-existence of finite order transcendental	entire solutions of the certain non-linear differential-difference equations. We also  study about conjectures posed by Rong et al. and Chen et al. \\  
\end{abstract}
	
\section{\textbf{Introduction and Main Results}}
In this paper, we use standard notations of Nevanlinna theory such as $T(r,f), N(r,f)$ and $m(r,f)$ to denote characteristic function, counting function and proximity function of $f$ respectively. We also use fundamental results of Nevanlinna theory, see \cite{ilaine, yanglo}. The terms $\rho(f), \rho_2(f)$ denote order of $f$ and hyper-order of $f$ respectively.\\
If $f$ is a meromorphic function, then the quantities  which are of growth $o(T(r,f))$ as $r\to\infty$, outside a set of finite linear measure, are denoted by $S(r,f)$. We say that a meromorphic function $g(z)$ is a small function of $f(z)$ if $T(r,g)=S(r,f)$ and vice versa. We also note that the finite sum of quantities of type $S(r,f)$ is again of type $S(r,f)$. \\
It is quite interesting to check existence or non-existence of solutions of non-linear differential equations or difference equations or differential-difference equations. Many researchers are working in this field, see \cite{li, ly, llsy, luwu, chw, zl}.\\
Recently Rong and Xu \cite{rx} studied a problem in this direction which is as follows:	
\begin{theorem}\label{baseth}
Suppose that $\alpha_1, \alpha_2, p_1, p_2$ and $\beta$ be non zero constants satisfying $\alpha_1\neq\alpha_2$. Suppose that $k\geq 0$ and $n\geq2$ are integers and $a(z)$ is a non zero polynomial. If $f(z)$ is a transcendental entire solution with $\rho_2(f)<1$ of the differential-difference equation
$$f^n(z)+a(z)f^{(k)}(z+\beta)=p_1e^{\alpha_{1}z}+p_2e^{\alpha_{2}z},$$
then we have $\rho(f)=1$, $a(z)$ must be a constant and one of the following relations holds:
\begin{enumerate}[$(i)$]
\item $f(z)=c_1e^{\frac{\alpha_1z}{n}},~   ac_1({\frac{\alpha_1}{n}})^ke^{\frac{\alpha_1\beta}{n}}=p_2, ~ \alpha_1=n\alpha_2$,
	
\item $f(z)=c_2e^{\frac{\alpha_2z}{n}},~   ac_2({\frac{\alpha_2}{n}})^ke^{\frac{\alpha_2\beta}{n}}=p_1, ~ \alpha_2=n\alpha_1$,\\
where $c_1, c_2$ are constants satisfying $c_i^n=p_i$, ~ $i=1,2.$
\item If $n=2$, we have $T(r,f)\leq N_{1)}
(r,1/f)+T(r,\psi)+S(r,f)$, where $N_{1)}(r,1/f)$ denotes the counting function corresponding
to simple zeros of $f$ and $(0\not\equiv)\psi=\alpha_1\alpha_2f^2-n(\alpha_1+\alpha_2)ff'+nff''+n(n-1)(f')^2$. If $n=3$, we have $T(r,f)=N_{1)}(r,1/f)+S(r,f)$. If $n\geq4$, we only have the cases (i) and (ii).
\end{enumerate}
\end{theorem}

In the above theorem, Rong and Xu \cite{rx} conjectured that case (iii) can be removed if $N(r,1/f)=S(r,f)$ for $n=2,3.$	\\
We give affirmative answer for $n=2$ in the following Theorem. 
\begin{theorem}\label{mainth1}
Suppose that $p_{1}, p_{2}, \alpha_{1},  \alpha_{2}$ are non zero constants such that $\alpha_{1}\neq\alpha_{2}$, $\beta$ is any constant, $a(z)$ is a non zero polynomial and $k\geq 0$ is an integer. If $f(z)$ is a transcendental entire solution with $\rho_2(f)<1$ of the differential-difference equation
\begin{equation}\label{mdde}
f^2(z)+a(z)f^{(k)}(z+\beta)=p_1e^{\alpha_{1}z}+p_2e^{\alpha_{2}z}
\end{equation}
and satisfies $N(r,\frac{1}{f})=S(r,f)$, then $\rho(f)=1$, $a(z)$ must be a constant and any one of the following holds:
\begin{enumerate}[(i)]
\item $f(z)=c_1e^{\frac{\alpha_1z}{2}},~   ac_1({\frac{\alpha_1}{2}})^ke^{\frac{\alpha_1\beta}{2}}=p_2, ~ \alpha_1=2\alpha_2$,

\item $f(z)=c_2e^{\frac{\alpha_2z}{2}},~   ac_2({\frac{\alpha_2}{2}})^ke^{\frac{\alpha_2\beta}{2}}=p_1, ~ \alpha_2=2\alpha_1$,

\end{enumerate}
where $c_1, c_2$ are constants satisfying $c_i^2=p_i,i=1,2$.
\end{theorem}
As we know that an entire function with finite order of growth has zero hyper-order of growth. Therefore the following corollaries are immediate consequences of the above theorem. 
\begin{corollary}
In the hypothesis of Theorem \ref{mainth1}, if $f$ is a finite order of transcendental entire solution of \eqref{mdde}, satisfying $N(r,1/f)=S(r,f)$. Then the same conclusions hold.
\end{corollary}
\begin{corollary}
In the hypothesis of Theorem \ref{mainth1}, if $f$ is a finite order regular growth of transcendental entire solution of \eqref{mdde}, having $0$ as a Picard exceptional value. Then the same conclusions hold.
\end{corollary}

\begin{example}
The function $f(z)=3e^{2z}$ satisfies
$$f^2(z)+\frac{1}{2}f''(z+\log2)=9e^{4z}+24e^{2z},$$
where $f(z)$ holds conclusion $(i)$ of Theorem \ref{mainth1}.\\
\end{example}

By some further analysis, we prove a result which shows that under certain conditions, equation \eqref{maineq2} does not have finite order of transcendental entire solution.

\begin{theorem}\label{mainth3}
Let us consider a particular kind of differential-difference equation
\begin{equation}\label{maineq2}
f^2(z)+a(z)f^{(k)}(z+\beta)=\alpha(z),
\end{equation}
where $a(z)$ is a small function of $f(z)$, $k\geq 0$ is an integer, $\beta$ is a constant and $\alpha(z)$ is a meromorphic function such that
$T(r,\alpha)=\gamma T(r,f)+S(r,f)$, where $0\leq\gamma<1$. Then there does not exist a finite order transcendental entire solution of \eqref{maineq2}. \\
\end{theorem}	

Recently Chen et al. \cite{chw} investigated the entire solutions with finite order of the following differential-difference equation
\begin{equation}\label{chweq}
f^n(z)+wf^{n-1}(z)f'(z)+q(z)e^{Q(z)}f(z+c)=p_1e^{\lambda z}+p_2e^{-\lambda z},
\end{equation}
where $n$ is an integer, $w$ is a constant, $p_1, p_2$ and $\lambda$ are non-zero constants, $q\not\equiv 0$ is a polynomial and $Q$ is a non constant polynomial. They studied the following theorem:
\begin{theorem}\label{conjectureth}
Suppose that $f$ is a finite order transcendental entire solution of \eqref{chweq}, then the following conclusions hold:
\begin{enumerate}[$(i)$]
\item If $n\geq 4$ for $w\neq0$ and $n\geq3$ for $w=0$, then every solution $f$ satisfies
$\rho(f)=\deg(Q)=1$.
\item If $n\geq 1$ and $f$ is a solution of \eqref{chweq}, which belongs to $\Gamma_0=\{e^{\gamma(z)}:\gamma(z)$ is a non-constant polynomial\}, then
$$f(z)=e^{\lambda z/n+a}, \qquad \qquad  Q(z)=-\frac{n+1}{n}\lambda z+b$$
or
$$f(z)=e^{-\lambda z/n+a}, \qquad \qquad  Q(z)=\frac{n+1}{n}\lambda z+b,$$
where $a, b\in \mathbb{C}$.
\end{enumerate}
\end{theorem}
In the Theorem \ref{conjectureth}, Chen et al. \cite{chw} conjectured that conclusion $(i)$ will also hold for $n=2$. In the next theorem, we are giving partial answer to their conjecture by adding $N(r,1/f)=S(r,f)$ condition in the hypothesis. 

\begin{theorem}\label{mainth4}
If $f$ is a transcendental entire solution of finite order of
\begin{equation}\label{edde}
f^2(z)+wf(z)f'(z)+q(z)e^{Q(z)}f(z+c)=p_1e^{\lambda z}+p_2e^{-\lambda z}
\end{equation}
and satisfying $N(r,1/f)=S(r,f)$, then every solution $f$ satisfies $\rho(f)=\deg(Q)=1$.
\end{theorem}

\begin{example}
The function $f(z)=e^z$ with $N(r,1/f)=S(r,f)$ satisfies 
$$f^2+2ff'+\frac{1}{5} e^{-3z}f(z+\log 10)=3e^{2z}+2e^{-2z},$$
here $\rho(f)=\deg(Q)=1$.\\
\end{example}

In 2014, Liu et al.\cite{llsy} studied a problem based on certain type of non-linear difference equation which is as follows:
\begin{theorem}
Suppose that $p_1, p_2, \alpha_1, \alpha_2$ be non-zero constants,  $q(z)$ is a polynomial and $n\geq4$ is an integer. If $f$ is a finite order entire function to the following equation
\begin{equation}\label{nde}
f^n+q(z)\Delta f(z)=p_1e^{\alpha_1z}+p_2e^{\alpha_2z},
\end{equation}
then $q$ is a constant, and one of the following relations holds:
\begin{enumerate}[$(i)$]
\item  $f(z)=c_1e^{\frac{\alpha_1z}{n}},~   c_1(e^{\frac{\alpha_1}{n}}-1)q=p_2, ~ \alpha_1=n\alpha_2$
\item $f(z)=c_1e^{\frac{\alpha_2z}{n}},~   c_1(e^{\frac{\alpha_2}{n}}-1)q=p_1, ~ \alpha_2=n\alpha_1$,
\end{enumerate}
where $c_1, c_2$ are constants satisfying $c_1^n= p_1, c_2^n = p_2$.
\end{theorem}
Later on, Zhang et al.\cite{zlly} studied equation \eqref{nde} for $n=3$. Further, L$\ddot{u}$ et al.\cite{luwu} studied a non-linear differential-difference equation for $n=3$, which was slight general form of \eqref{nde}. Recently Xu and Rong studied a non-linear difference equation for $n=2$, see \cite{xr}.

%\begin{equation}\label{3degeq1}
%f^3+L_1(z,f)=p_1e^{\alpha_1z}+p_2e^{\alpha_2z},
%\end{equation} 

%where $p_1, p_2, \alpha_1, \alpha_2$ be non-zero constants such that $\alpha_1\neq \alpha_2$, $L_1(z,f)$ denote a differential-difference polynomial in $f$ of degree one with small functions of $f$ as its coefficients. Note that \eqref{3degeq1} is slight general form of equation \eqref{nde}.\\
%\begin{theorem}
%Suppose that and $p_1, p_2, \alpha_1, \alpha_2$ be non-zero constants such that $\alpha_1\neq \alpha_2$. If $f$ is an entire solution with $\rho_2(f)<1$ to the following equation
%$$f^3+L_1(z,f)=p_1e^{\alpha_1z}+p_2e^{\alpha_2z},$$
%then one of the following relations holds
%\begin{enumerate}[$(i)$]
%\item $f(z)=c_1\exp(
%\alpha_1z/3)+c_2\exp(\alpha_2z/3)$, $c_i^3=p_i; i=1,2$ and $\alpha_1+\alpha_2=0$, where $c_1, c_2$ are non-zero constants
%\item $f^3(z)=(p_1-c_1)\exp(\alpha_1z)$ and $L_1(z,f)=c_1\exp(\alpha_1z)+p_2\exp(\alpha_2z)$, where $c_1$ is a constant
%\item $f^3(z)=(p_2-c_2)\exp(\alpha_2z)$ and $L_1(z,f)=p_1\exp(\alpha_1z)+c_2\exp(\alpha_2z)$, where $c_2$ is a constant.
%\end{enumerate}
%\end{theorem}

Motivated by the above results, we are dealing with non-linear difference equation, which is as follows:
\begin{theorem}\label{mainth2}
Suppose that $p_{1}, p_{2}, \alpha_{1},  \alpha_{2}$ are non zero constants such that $\alpha_{1}\neq\alpha_{2}$, $L(z,f)$ is a linear difference polynomial with constant coefficients. If $f(z)$ is a transcendental entire solution with $\rho_2(f)<1$ of the difference equation
\begin{equation}\label{mde}
f^2(z)+L(z,f)=p_1e^{\alpha_{1}z}+p_2e^{\alpha_{2}z},
\end{equation}
then $\rho(f)=1$ and any one of the following holds:
\begin{enumerate}[(i)]
\item $T(r,f)\leq 2N\left(r,\frac{1}{f}\right)+S(r,f)$
\item $f(z)=c_1e^{\alpha_1z/2},~ c_1^2=p_1$
\item $f(z)=c_2e^{\alpha_2z/2},~ c_2^2=p_2$,  
\end{enumerate}
where $c_1, c_2$ are constants.
\end{theorem}

\begin{remark}
In the above hypothesis, if we add $N(r,1/f)=S(r,f)$, then case $(i)$ will be removed.
\end{remark}

\begin{example}
The function $f(z)=2e^{z}$ satisfies the difference equation
$$f^2(z)+3f(z+2\pi i)=4e^{2z}+6e^{z},$$
which shows that $f(z)$ holds conclusion $(ii)$ of Theorem \ref{mainth2}.\\
\end{example}

\begin{example}
The function $f(z)=-8+4e^{\pi iz}+8e^{-\pi iz}$ satisfies 
$$f^2(z)+3f(z+2)+13f(z+8)=16e^{2\pi iz}+64e^{-2\pi iz},$$
which shows that $f(z)$ holds conclusion $(i)$ of Theorem \ref{mainth2}.\\
\end{example}

%\begin{example}	
%The function $f(z)=-6+3\sqrt{2}e^{\pi iz}+3\sqrt{2}e^{-\pi iz}$ satisfies
%$$f^2(z)+6f(z+2)+6f(z+4)=18(e^{2\pi iz}+e^{-2\pi iz}),$$
%which also shows that $f(z)$ holds conclusion $(i)$ of Theorem \ref{mainth2}.\\
%\end{example}

\begin{example}
The function $f(z)=5e^{3z/2}$ satisfies
$$f^2(z)+2f(z+2)+5f(z+3)=p_1e^{3z/2}+25e^{3z},$$
where $p_1=10e^3+25e^{9/2}$. Here $f(z)$ holds conclusion $(iii)$ of Theorem \ref{mainth2}.\\
\end{example}

\section{\textbf{Preliminary Results}}

This section includes those results which will be used to prove our main theorems.

\begin{lemma}\rm{\cite{ilaine}} \label{il}
Suppose that $f$ is a transcendental meromorphic function and $k\geq1$ be an integer. Then
$$m\left(r,\frac{f^{(k)}}{f}\right)=S(r,f),$$
and if $f$ is a finite order transcendental meromorphic function, then
$$m\left(r,\frac{f^{(k)}}{f}\right)=O(\log r).$$\\
\end{lemma}

%\begin{lemma}
%Suppose that $c_1, c_2$ are complex numbers such that $c_1\neq n_2$ and $f$ is a meromorphic function with $\rho(f)<\infty$. Then for each $\epsilon>0$, we have
%$$m\left(r,\frac{f(z+c_1)}{f(z+c_2)}\right)=O(r^{\rho-1+\epsilon})$$
%\end{lemma}

The following lemma is a difference analogue of the lemma on the logarithmic derivative for finite order meromorphic functions. It was proved by Halburd and Korhonen \cite{hk1,hk2} and Chiang and Feng \cite{cf} independently.
\begin{lemma}\rm{\cite{cf}}\label{hk}
Suppose that $f$ is a meromorphic function with $\rho(f)<\infty$ and $c_1,c_2\in\mathbb{C}$ such that $c_1\neq c_2$, then  for each $\epsilon>0$, we have
$$m\left(r,\frac{f(z+c_1)}{f(z+c_2)}\right) =O(r^{\rho-1+\epsilon}).$$
\end{lemma}

\begin{remark} \rm{\cite{hk1,hk2}}
Suppose that $f$ is a non-constant meromorphic function with $\rho(f)<\infty$ and $c\in\mathbb{C}$, then
$$m\left(r,\frac{f(z+c)}{f(z)}\right)=S(r,f),$$
outside a possible exceptional set of finite logarithmic measure.\\
The conclusion of above remark also holds when $f$ is a meromorphic function of hyper-order $\rho_2(f)<1$. See \cite{hkt}.
\end{remark}
%\begin{lemma}\rm{\cite{ilaine}}
%Suppose that $f$ is a transcendental meromorphic solution of 
%$$f^nP(z,f)=Q(z,f),$$
%where $P(z,f)$ and $Q(z,f)$ are polynomials in $f(z)$ and its derivatives with meromorphic
%coefficients, say $\{a_\lambda : \lambda\in I\}$, satisfying $m(r,a_\lambda)=S(r, f)$ for all $\lambda\in I$. If the total degree of
%$Q(z,f)$ as a polynomial in $f(z)$ and its derivatives is at most $n$, then
%$$m(r,P(z,f))=S(r,f).$$
%\end{lemma}

%\begin{lemma}\label{differencelemma}
%Suppose that $f(z)$ is a non-constant meromorphic solution of finite order of 
%$$ f^nP(z,f)=Q(z,f),$$
%where $P(z,f)$ and $Q(z,f)$ are polynomials in $f(z)$, its derivatives and its shifts $f(z+c)$ with small meromorphic coefficients. If the total degree of $Q(z,f)$ is at most $n$, then
%$$m(r,P(z,f))=S(r,f)$$
%for all $r$ outside a possible set of finite logarithmic measure.
%\end{lemma}
Next lemma plays an important role to prove results in differential or differential-difference equations which is as follows: 
\begin{lemma}\label{imp}\rm\cite{yybook}
let $f_1,f_2,...,f_n (n\geq2)$ be meromorphic
functions and $g_1,g_2,...,g_n$ be entire functions satisfying the following conditions:
\begin{enumerate}
\item $\sum_{i=1}^{n}f_ie^{g_i}\equiv 0$.
\item $g_j-g_k$ are not constants for $1\leq j<k\leq n$.
\item For $1\leq i\leq n, 1\leq h<k\leq n$,\\
$T(r, f_i)=o(T(r,e^{(g_h-g_k)}))$ as $r\to\infty$, outside a set of finite logarithmic measure.
\end{enumerate}
 Then $f_i\equiv 0$, $(i=1,2,...,n).$

\end{lemma}

%\begin{lemma}\label{cg}
%Suppose that $R(z)$ is a non zero polynomial and $\gamma$ is a non zero constant. Then the differential equation
%$$\gamma f(z)-2f'(z)=H(z)$$ 
%has a special solution $c_0(z)$ which is a non zero polynomial.
%\end{lemma}

\section{\textbf{Proof of Theorems}}
\begin{proof}[\textbf{\underline{Proof of Theorem \ref{mainth1}}}]
Suppose that $f(z)$ is a transcendental entire solution with $\rho_2(f)<1$. By the same method as in [\cite{rx},Theorem 7], we have $\rho(f)=1$.\\
Let $P=a(z)f^{(k)}(z+\beta)$ and differentiating \eqref{mdde}, we have
\begin{equation}\label{dmdde}
2ff'+P'=P_1\alpha_1e^{\alpha_1z}+P_2\alpha_2e^{\alpha_2z}.
\end{equation}
Eliminating $e^{\alpha_1z}$ and $e^{\alpha_2z}$ by using \eqref{mdde} and \eqref{dmdde}, respectively, we obtained
\begin{equation}\label{e2}
\alpha_1f^2+\alpha_1P-2ff'-P'=(\alpha_1-\alpha_2)p_2e^{\alpha_2z},
\end{equation}
\begin{equation}\label{e1}
\alpha_2f^2+\alpha_2P-2ff'-P'=(\alpha_2-\alpha_1)p_2e^{\alpha_1z}.
\end{equation}
Now differentiating \eqref{e2}, we have
\begin{equation}\label{de2}
2\alpha_1ff'+\alpha_1P'-2(f')^2-2ff''-P''=(\alpha_1-\alpha_2)p_2\alpha_2e^{\alpha_2z}.
\end{equation}
From \eqref{e2} and \eqref{de2}, we have
	
\begin{equation}\label{keyeqq}
\psi=Q,
\end{equation}
where
\begin{equation}\label{keyeqq1}
\psi=\alpha_1\alpha_2f^2-2(\alpha_1+\alpha_2)ff'+2ff''+2(f')^2
\end{equation}
and 
\begin{equation}\label{keyeqq2}
Q=-\alpha_1\alpha_2P+(\alpha_1+\alpha_2)P'-P''
\end{equation}
Now we discuss two cases:
\begin{enumerate}
\item If $\psi\equiv 0$, proceeding on similar lines to Theorem 7 in \cite{rx}, we have desired conclusions.\\

\item If $\psi\not\equiv 0$, using \eqref{keyeqq}-\eqref{keyeqq2}, Lemma \ref{il}, and Lemma \ref{hk}, we obtain
\begin{equation}\label{mrfeq}
m\left(r,\frac{\psi}{f}\right)=m\left(r,\frac{Q}{f}\right)=S(r,f), \qquad \qquad m\left(r,\frac{\psi}{f^2}\right)=S(r,f).
\end{equation}
%Given that $N\left(r,\frac{1}{f}\right)=S(r,f)$ and using \eqref{mrfeq}, we obtain\\
%\begin{align*}
%2T(r,f)=T(r,f^2)&=m\left(r,\frac{1}{f^2}\right) +S(r,f) \\
%&\leq m\left(r,\frac{\psi}{f^2}\right)+ m\left( r,\frac{1}{\psi}\right) +S(r,f) \\
%&\leq m(r,\psi)+S(r,f) \\
%&\leq m\left( r,\frac{\psi}{f}\right) +m(r,f)+S(r,f) \\
%&\leq m(r,f)+S(r,f)=T(r,f)+S(r,f)
%\end{align*}
%This yields to $T(r,f)\leq S(r,f)$, which is a contradiction. This also proves that Rong and Xu's conjecture was true for $n=2$.\\

%OR\\

Also from \eqref{keyeqq1}, we have
\begin{equation}
\frac{1}{f^2}=\frac{1}{\psi} \left[-\alpha_1\alpha_2-2(\alpha_1+\alpha_2)\frac{f'}{f}+2\frac{f''}{f}+2\left(\frac{f'}{f}\right)^2\right].
\end{equation} 
Using Lemma \ref{il}, we have
\begin{equation}\label{1-eq}
2m\left(r,\frac{1}{f}\right)=m\left(r,\frac{1}{f^2}\right)\leq m\left(r,\frac{1}{\psi}\right)+S(r,f).
\end{equation}
%Given that $N\left(r,\frac{1}{f}\right)=S(r,f)$, using \eqref{1-eq}, \eqref{mrfeq} and first fundamental theorem of Nevanlinna, we have
%\begin{align*}
%2T(r,f)=T(r.f^2)&=m\left(r,\frac{1}{f^2}\right)+N\left(r,\frac{1}{f^2}\right)+O(1)\\
%&\leq m\left(r,\frac{1}{\psi}\right)+S(r,f)\\
%&\leq T(r,\psi)+S(r,f)=m(r,\psi)+S(r,f)\\
%&\leq m\left(r,\frac{\psi}{f}\right)+m(r,f)+S(r,f)\\
%&\leq m(r,f)+S(r,f)=T(r,f)+S(r,f)
%\end{align*}
%This implies $T(r,f)\leq S(r,f)$, which is a contradiction. Now we also prove that Rong and Xu's conjecture was true for $n=2$.\\
Given that $N\left(r,\frac{1}{f}\right)=S(r,f)$, using \eqref{1-eq} and first fundamental theorem of Nevanlinna, we have  
\begin{align*}
T(r,f)&=m\left(r,\frac{1}{f}\right)+N\left(r,\frac{1}{f}\right)+O(1)\\
&\leq 2m\left(r,\frac{1}{f}\right)+S(r,f)\\
&\leq m\left(r,\frac{1}{\psi}\right)+S(r,f)\leq T(r,\psi)+S(r,f).
\end{align*}
This implies
\begin{equation}\label{conjeceq}
T(r,f)\leq T(r,\psi)+S(r,f).
\end{equation}
Using \eqref{mrfeq}, we have
\begin{align*}
T(r,\psi)=m(r,\psi) &\leq m(r,\frac{\psi}{f})+m(r,f)\\
&=S(r,f)+m(r,f)=T(r,f)+S(r,f).
\end{align*}
Thus
\begin{equation}\label{conjeckeyeq}
T(r,\psi)\leq T(r,f)+S(r,f).
\end{equation} 
From the above two inequalities \eqref{conjeceq} and \eqref{conjeckeyeq}, we have
$$T(r,\psi)=T(r,f)+S(r,f).$$ 
Now using \eqref{mrfeq}, we have
\begin{align*}
2T(r,\psi)&=2T(r,f)+S(r,f)\\
&=m(r,\frac{1}{f^2})+S(r,f) \qquad \qquad \left(Since, N\left(r,\frac{1}{f}\right)=S(r,f)\right)\\
&\leq m(r,\frac{\psi}{f^2})+m(r,\frac{1}{\psi})+S(r,f)\\
&\leq T(r,\psi)+S(r,f).
\end{align*}
This implies $T(r,\psi)\leq S(r,f)$. From \eqref{conjeceq}, we obtain $T(r,f)\leq S(r,f)$, which is a contradiction. This also proves that Rong and Xu's conjecture was true for $n=2$.
\end{enumerate}
\end{proof}

\begin{proof}[\textbf{\underline{Proof of Theorem \ref{mainth3}}}]
We prove this theorem by contradiction. Suppose that $f$ is a transcendental entire solution of \eqref{maineq2}, then using Lemma \ref{il}, Lemma \ref{hk} and $T(r,\alpha(z))=\gamma T(r,f)+S(r,f); 0\leq\gamma<1$, we have
\begin{align*}
2T(r,f)=T(r,f^2)&=T(r,\alpha(z)-a(z)f^{(k)}(z+\beta))\\
&\leq T(r,\alpha(z))+T(r,a(z))+T(r,f^{(k)}(z+\beta))\\
&\leq \gamma T(r,f)+S(r,f)+T\left(r,\frac{f^{(k)}(z+\beta)}{f(z)}\right)+T(r,f)\\
&\leq \gamma T(r,f)+S(r,f)+T\left(r,\frac{f^{(k)}(z+\beta)}{f(z+\beta)}\right)+T\left(r,\frac{f(z+\beta)}{f(z)}\right)+T(r,f)\\
&\leq \gamma T(r,f)+T(r,f)+S(r,f).
\end{align*} 
This implies
$$T(r,f)\leq \gamma T(r,f)+S(r,f),$$
which is a contradiction. This completes the proof.\\
\end{proof}

\begin{proof}[\textbf{\underline{Proof of Theorem \ref{mainth4}}}]
Suppose that $f$ is a transcendental entire solution of finite order of \eqref{edde}. To prove the desired conclusion, we discuss following cases:\\
	
\underline{\textbf{Case I:}} If $\rho(f)<1$, this implies $T(r,f)=S(r,e^{\lambda z})$. From 
\eqref{edde}, Lemma \ref{il} and Lemma \ref{hk}, we have
\begin{align*}
T(r,e^{Q(z)})=m(r,e^{Q(z)})&=m\left(r,\frac{p_1e^{\lambda z}+p_2e^{-\lambda z}-wff'-f^2}{q(z)f(z+c)}\right)\\
&\leq m\left(r,\frac{1}{q(z)f(z+c)}\right)+m(r,p_1e^{\lambda z}+p_2e^{-\lambda z})+m(r,wff')+m(r,f^2)\\
&\leq m\left(r,\frac{f(z)}{q(z)f(z+c)}\right)+m\left(r,\frac{1}{f(z)}\right)+2m(r,e^{\lambda z})+S(r,e^{\lambda z})\\
&\leq 2m(r,e^{\lambda z})+S(r,e^{\lambda z})=2T(r,e^{\lambda z})+S(r,e^{\lambda z}).
\end{align*}
This implies $\deg(Q)=\rho(e^{Q(z)})\leq 1$. Also we know that $\deg(Q)\geq 1$, thus $\deg(Q)=1$.\\
Let us consider $Q(z)=az+b$, where $a\in\mathbb{C}\setminus\{0\}$ and $b\in\mathbb{C}$.\\
Rewrite \eqref{edde} as 
\begin{equation}\label{4-1eq}
f^2+wff'+q(z)e^{(az+b)}f(z+c)=p_1e^{\lambda z}+p_2e^{-\lambda z}.
\end{equation}
Differentiating \eqref{4-1eq}, we have
\begin{equation}\label{4-2eq}
2ff'+wff''+w(f')^2+Le^{(az+b)}=\lambda(p_1e^{\lambda z}-p_2e^{-\lambda z}),
\end{equation}
where $L=q'(z)f(z+c)+aq(z)f(z+c)+q(z)f'(z+c)$. For convenience, further we use $f_c$ and $q$ in place of $f(z+c)$ and $q(z)$ respectively.  \\
	
Eliminating $e^{\lambda z}, e^{-\lambda z}$ from \eqref{4-1eq} and \eqref{4-2eq}, we have
\begin{equation}\label{4-3eq}
A(z)e^{2(az+b)}+B(z)e^{(az+b)}+C(z)\equiv 0,
\end{equation}
where
\begin{align*}
A(z)&=\lambda^2q^2f^2_c-L^2,\\
B(z)&=2\lambda^2(qf^2f_c+wqff'f_c)-2A(2ff'+wf'^2+wff''),\\
C(z)&=\lambda^2(f^4+w^2f^2f'^2+2wf^3f')-(4f^2f'^2+w^2f'^4+w^2f^2f''^2+4wff'^3+\\
& \qquad \qquad 4wf^2f'f''+2w^2ff'^2f''+4\lambda^2p_1p_2).
\end{align*}
From \eqref{4-3eq} and Lemma \ref{imp}, we have
$$A(z)\equiv B(z)\equiv C(z)\equiv 0.$$
It follows from $C(z)\equiv 0$ that  
\begin{equation}\label{4-keyeq1}
\phi(f)= 4\lambda^2p_1p_2,
\end{equation}
where
\begin{equation}\label{4-4eq}
\phi(f)=\lambda^2(f^4+w^2f^2f'^2+2wf^3f')-(4f^2f'^2+w^2f'^4+w^2f^2f''^2+4wff'^3+4wf^2f'f''+2w^2ff'^2f'').
\end{equation}
From \eqref{4-keyeq1}, \eqref{4-4eq}  and Lemma \ref{il}, we obtain
\begin{equation}\label{4-keyeq2}
m(r,\phi)=S(r,f), \qquad \qquad m\left(r,\frac{\phi}{f^4}\right)=S(r,f).
\end{equation}   
Since $N(r,1/f)=S(r,f)$, using \eqref{4-keyeq2} and first fundamental theorem of Nevanlinna, we have
\begin{align*}
4T(r,f)=T(r,f^4)&=m(r,\frac{1}{f^4})+S(r,f)\\
&\leq  m(r,\frac{\phi}{f^4})+m(r,\frac{1}{\phi})+S(r,f)\\
&\leq m(r,\phi)+S(r,f)=S(r,f).
\end{align*}
This implies $T(r,f)\leq S(r,f)$, which is a contradiction.\\
	
\underline{\textbf{Case II:}} If $\rho(f)>1$, since $\rho(p_1e^{\lambda z}+p_2e^{-\lambda z})=1$, then $T(r,p_1e^{\lambda z}+p_2e^{-\lambda z})=S(r,f)$.\
	
Differentiating \eqref{edde}, we have
\begin{equation}\label{4-5eq}
2ff'+wff''+wf'^2+Re^{Q(z)}=\lambda(p_1e^{\lambda z}-p_2e^{-\lambda z}),
\end{equation}

where $R=q'f_c+qQ'f_c+qf'_c$.\\

From \eqref{4-5eq}, substituting the value of $e^{Q(z)}$ into \eqref{edde} gives
\begin{equation}\label{4-keyeq3}
\psi(f)=P,
\end{equation}

where
\begin{equation}\label{4-6eq}
\psi(f)=Rf^2+wRff'-qf_c(2ff'+wf'^2+wff'')
\end{equation}

and 
\begin{equation}\label{4-7eq}
P=(R-\lambda qf_c)e^{\lambda z}+(R+\lambda qf_c)e^{-\lambda z}.
\end{equation}

Since $\rho(f)>1$, this implies $T(r,e^{\lambda z})=S(r,f)$ and $T(r,e^{-\lambda z})=S(r,f)$.\\

Using \eqref{4-keyeq3}, Lemma \ref{il} and Lemma \ref{hk}, we get
\begin{equation}\label{4-keyeq4}
m\left(r,\frac{\psi}{f}\right)=m\left(r,\frac{P}{f}\right)=S(r,f), \qquad \qquad m\left(r,\frac{\psi}{f^3}\right)=S(r,f).
\end{equation}

Since $N(r,1/f)=S(r,f)$, using \eqref{4-keyeq4} and first fundamental theorem of Nevanlinna, we have
\begin{align*}
3T(r,f)=T(r,f^3)&=m(r,\frac{1}{f^3})+S(r,f)\\
&\leq m(r,\frac{\psi}{f^3})+m(r,\frac{1}{\psi})+S(r,f)\\
&\leq T(r,\psi)+S(r,f)=m(r,\psi)+S(r,f)\\
&\leq m(r,\frac{\psi}{f})+m(r,f)+S(r,f)\\	&=m(r,f)+S(r,f)=T(r,f)+S(r,f).
\end{align*}

This implies $T(r,f)\leq S(r,f)$, which is a contradiction.\\
	
\underline{\textbf{Case III:}} If $\rho(f)=1$, then from \eqref{edde}, Lemma \ref{il} and Lemma \ref{hk}, we have
\begin{align*}
T(r,e^{Q(z)})&=m(r,e^{Q(z)})\\
&=m\left(r,\frac{p_1e^{\lambda z}+p_2e^{-\lambda z}-f^2-wff'}{q(z)f(z+c)}\right)\\
&\leq m\left(r,\frac{1}{q(z)f(z+c)}\right)+2m(r,e^{\lambda z})+m(r,f(f+wf'))+S(r,f)\\
&\leq m\left(r,\frac{f(z)}{q(z)f(z+c)}\right)+m\left(r,\frac{1}{f(z)}\right)+2m(r,e^{\lambda z})+m\left(r,\frac{f+wf'}{f}\right)+\\
 & \qquad \qquad \qquad \qquad 2m(r,f)+S(r,f)\\
&\leq 2m(r,e^{\lambda z})+3m(r,f)+S(r,f)=2T(r,e^{\lambda z})+3T(r,f)+S(r,f).
\end{align*}
This implies $\deg(Q)=\rho(e^{Q(z)})\leq \max\{\rho( e^{\lambda z}),\rho(f)\}=1$. Also we know that $\deg(Q)\geq 1$, thus we have $\deg(Q)=1$.\\
\end{proof}

\begin{proof}[\textbf{\underline{Proof of Theorem \ref{mainth2}}}]
Suppose that $f$ is a transcendental entire solution with $\rho_2(f)<1$. For convenience, we use $L$ in place of $L(r,f)$. From \eqref{mde} and Lemma \ref{hk}, we have
\begin{align*}
T(r,p_1e^{\alpha_{1}z}+p_2e^{\alpha_{2}z})=T(r,f^2+L)&\leq T(r,f^2)+T(r,L)\\
 &=m(r,f^2)+m(r,L)\\
&\leq m(r,f^2)+m\left(r,\frac{L}{f}\right)+m(r,f)\\
&\leq 3m(r,f)+S(r,f)=3T(r,f)+S(r,f)
\end{align*}

and
\begin{align*}
T(r,p_1e^{\alpha_{1}z}+p_2e^{\alpha_{2}z})&=T(r,f^2+L)\\
&\geq m(r,f^2)-m(r,L)\\
&=2m(r,f)-m\left(r,\frac{L}{f}f\right)\\
&\geq 2m(r,f)-m(r,f)+S(r,f)\\
&=m(r,f)+S(r,f)=T(r,f)+S(r,f).
\end{align*}

From the above two inequalities, we get
\begin{equation*}
T(r,f)+S(r,f)\leq T(r,p_1e^{\alpha_{1}z}+p_2e^{\alpha_{2}z})\leq 3T(r,f)+S(r,f).
\end{equation*} 

Since $\rho(p_1e^{\alpha_{1}z}+p_2e^{\alpha_{2}z})=1$, this implies $\rho(f)=1$.\\

Differentiating \eqref{mde}, we get
\begin{equation}\label{dmde}
2ff'+L'=p_1\alpha_1e^{\alpha_1z}+p_2\alpha_2e^{\alpha_2z}.
\end{equation}

Eliminating $e^{\alpha_1z}$ from \eqref{mde} and \eqref{dmde}, we obtain
\begin{equation}\label{eqe2}
\alpha_1f^2+\alpha_1L-2ff'-L'=(\alpha_1-\alpha_2)p_2e^{\alpha_2z}.
\end{equation}
%\begin{equation}\label{e1}
%\alpha_2f^2+\alpha_2P-2ff'-P'=(\alpha_2-\alpha_1)p_2e^{\alpha_1z}.
%\end{equation}

Now differentiating \eqref{eqe2}, we get
\begin{equation}\label{eqde2}
2\alpha_1ff'+\alpha_1L'-2(f')^2-2ff''-L''=(\alpha_1-\alpha_2)p_2\alpha_2e^{\alpha_2z}.
\end{equation}

From \eqref{eqe2} and \eqref{eqde2}, we have
	
\begin{equation}\label{keyeqq3}
\psi=Q,
\end{equation}

where
\begin{equation}\label{psieq}
\psi=\alpha_1\alpha_2f^2-2(\alpha_1+\alpha_2)ff'+2ff''+2(f')^2
\end{equation}

and 
\begin{equation}\label{qeq}
Q=-\alpha_1\alpha_2L+(\alpha_1+\alpha_2)L'-L''
\end{equation}

If $\psi\not\equiv0$, using \eqref{keyeqq3}, Lemma \ref{il}, and Lemma \ref{hk}, we obtain
\begin{equation}\label{mrfeq3}
m\left(r,\frac{\psi}{f}\right)=m\left(r,\frac{Q}{f}\right)=S(r,f), \qquad \qquad m\left(r,\frac{\psi}{f^2}\right)=S(r,f).
\end{equation}

Using \eqref{mrfeq3} and first fundamental theorem of Nevanlinna, we obtain
\begin{align*}
2T(r,f)=T(r,f^2)&=m\left(r,\frac{1}{f^2}\right) +N\left(r,\frac{1}{f^2}\right)+O(1) \\
&\leq m\left(r,\frac{\psi}{f^2}\right)+ m\left( r,\frac{1}{\psi}\right) +2N\left( r,\frac{1}{f}\right)+O(1) \\
&\leq S(r,f)+m(r,\psi)+2N\left( r,\frac{1}{f}\right) \\
&\leq m\left( r,\frac{\psi}{f}\right) +m(r,f)+2N\left( r,\frac{1}{f}\right)+S(r,f) \\
&\leq S(r,f)+m(r,f)+2N\left( r,\frac{1}{f}\right) \\
&=T(r,f)+2N\left( r,\frac{1}{f}\right)+S(r,f) .
\end{align*}

This yields to
	$$T(r,f)\leq 2N\left( r,\frac{1}{f}\right) +S(r,f),$$
which is one of the conclusions.\\
	
If $\psi\equiv 0$, from \eqref{psieq}, we have
$$\alpha_1\alpha_2f^2-2(\alpha_1+\alpha_2)ff'+2ff''+2(f')^2=0.$$

Now dividing above equation by $f^2$ and substituting $\frac{f'}{f}=s$, we have
$$s'+2s^2-(\alpha_1+\alpha_2)s+\frac{\alpha_1\alpha_2}{2}=0.$$

This equation has two constant solutions $s_1=\frac{\alpha_1}{2}$ and $s_2=\frac{\alpha_2}{2}$.\\

Let $s\neq\frac{\alpha_1}{2}$ and $s\neq\frac{\alpha_2}{2}$, we have
$$\left(\frac{s'}{s-s_1}-\frac{s'}{s-s_2}\right)=-2(s_1-s_2).$$

Integrating the above equation, we have
$$\left(\frac{s-s_1}{s-s_2}\right) = \exp(-2(s_1-s_2)z+C), $$ 

where $C\in\mathbb{C}$.
This gives
$$s=\frac{s_2-s_1}{\exp(2(s_2-s_1)z+C)}+s_2=\frac{f'}{f}.$$

We observe that zeros of $f$ and  $\exp(2(s_2-s_1)z+C)$ are same. If $z_0$ is a zero of $f$ with multiplicity $n$, we have
$$n=Res\left(\frac{f'}{f},z_0\right)=Res\left(\frac{s_2-s_1}{\exp(2(s_2-s_1)z+C)}+s_2,z_0 \right)=\frac{1}{2}, $$

which is a contradiction.\\

Thus, if $s=\frac{\alpha_1}{2}$, we have $f=c_1e^{\frac{\alpha_1z}{2}}$, where $c_1\in\mathbb{C}$ such that $c_1^2=p_1$.\\

Similarly, if $s=\frac{\alpha_2}{2}$, we have $f=c_2e^{\frac{\alpha_2z}{2}}$, where $c_2\in\mathbb{C}$ such that $c_2^2=p_2$. \\

This completes the proof.
\end{proof}

\end{document}